\documentclass[12pt]{article}
\usepackage{amssymb}
\usepackage{graphicx} 
\newtheorem{lem}{Lemma} 
\newtheorem{thm}{Theorem}
\newtheorem{conj}{Conjecture}

\begin{document}

\title{Clique factors in powers of graphs}

\author{Ajit A. Diwan  and Aniruddha R. Joshi\\
Department of Computer Science and Engineering,\\
Indian Institute of Technology Bombay, Mumbai 400076, India.\\
email:\texttt{aad@cse.iitb.ac.in, aniruddhajoshi@cse.iitb.ac.in}}

\maketitle

\abstract{ The $k$th power of a graph $G$, denoted $G^k$, has the same vertex set as $G$, and two 
vertices are adjacent in $G^k$ if and only if there exists a path between them in $G$ of length at most 
$k$. A $K_r$-factor in a graph is a spanning subgraph in which every component is a complete
graph of order $r$. It is easy to show that for any connected graph $G$ of order divisible by
$r$, $G^{2r-2}$ contains a $K_r$-factor. This is best possible as there exist connected graphs $G$ 
of order divisible by $r$ such that $G^{2r-3}$ does not contain a $K_r$-factor. We conjecture that 
for any 2-connected graph $G$ of order divisible by $r$, $G^r$ contains a $K_r$-factor. This was known
for $r \le 3$ and  we prove it for $r = 4$. We prove a stronger statement that the vertex set of any 
2-connected graph $G$ of order $4k$ can be partitioned into $k$ parts of size $4$, such that the four 
vertices in any part are contained in a subtree of $G$ of order at most 5. More generally, we conjecture 
that for any partition of $n = n_1+n_2+\cdots+n_k$, the vertex set of any 2-connected graph $G$ of order 
$n$ can be partitioned into $k$ parts $V_1,V_2,\ldots,V_k$, such that $|V_i| = n_i$ and $V_i
\subseteq V(T_i)$ for some subtree $T_i$ of $G$ of order at most $n_i+1$, for $1 \le i \le k$.}

\section {Introduction}

One of the most well-studied topics in graph theory is to find the smallest value of a graph
parameter that ensures the existence of certain subgraphs. The parameters studied include
the number of edges, the minimum degree, the chromatic number and many others. Among the
subgraphs whose existence is studied are complete graphs of a given order, Hamiltonian cycles,
$k$ disjoint cycles and clique factors. Some of the classical theorems of this kind are 
Tur\'{a}n's theorem~\cite{T}, Dirac's theorem~\cite{D}, the Corr\'{a}di-Hajnal theorem~\cite{CH}
and the Hajnal-Szemer\'{e}di theorem~\cite{HS}. 

In this paper, we consider the smallest power of a graph that contains certain subgraphs. If
$G$ is a graph, then the $k$th power of $G$, denoted $G^k$, is the graph with the same vertex
set as $G$, and two vertices are adjacent in $G^k$ if and only if there exists a path
between them in $G$ of length at most $k$. If $G$ is a connected graph of order $n$, then
$G^{n-1}$ is the complete graph of order $n$, and thus for any graph $H$ of order at most $n$,
there exists a smallest value of $k$ such that $G^k$ contains $H$. It is well-known that for any 
connected graph $G$, $G^3$ contains a Hamiltonian cycle~\cite{K}~\cite{S} and if $G$ has even order, 
$G^2$ contains a perfect matching.  Fleischner's theorem~\cite{F} states that if $G$ is a 2-connected 
graph then $G^2$ is Hamiltonian. Several different proofs of stronger forms of this theorem 
have since been obtained, see for example~\cite{EF} and the references in it.

Let $r \ge 2$ be an integer. A $K_r$-factor in a graph $G$ is a spanning subgraph of $G$ in which 
every connected component is isomorphic to $K_r$, the complete graph on $r$ vertices. A trivial
necessary condition for a graph to contain a $K_r$-factor is that the number of vertices is
divisible by $r$ and we will only consider graphs satisfying this condition. If $G$ is a connected 
graph of order divisible by $r$, there exists a smallest integer $k$ such that $G^k$ contains a 
$K_r$-factor. It is not difficult to show that for any such graph $G$, $G^{2r-2}$ contains a 
$K_r$-factor and there exist connected graphs of order divisible by $r$ such that $G^{2r-3}$ does 
not contain a $K_r$-factor. We conjecture that if $G$ is a 2-connected graph of order divisible by 
$r$, then $G^r$ contains a $K_r$-factor. For $r=2$, this follows from the result for connected graphs, and 
it was proved for $r=3$ in~\cite{AD}. We prove it for $r=4$. This is the best possible, as there exist 
2-connected graphs $G$ of order divisible by $r$ such that $G^{r-1}$ does not contain a
$K_r$-factor.  

The proof is by induction, and we need to prove a stronger statement. A $K_4$-factor in $G^4$ is
equivalent to a partition of $V(G)$ into parts of size four such that any two vertices in the
same part are at distance at most four from each other in $G$. We impose a stronger condition on the
vertices in a part, and require that all four vertices in a part must be contained in a subtree
of $G$ with at most five vertices. This implies that for any two vertices in the same part, there
is a path of length at most four between them in the subtree that contains the part. We say a subset 
of vertices $A \subseteq V(G)$ is \emph{nearly connected} in a graph $G$ if there exists a subtree 
$T$ of $G$ of order at most $|A|+1$ such that $A \subseteq V(T)$. If $A$ is a nearly connected
subset of vertices in a graph $G$, then $A$ induces a complete subgraph in $G^{|A|}$.
We can now state our main result.

\begin{thm}
\label{main}
Let $G$ be any 2-connected graph of order $4k$ for some $k \ge 1$. Then $V(G)$ can be partitioned
into nearly connected subsets of size 4.
\end{thm}

We conjecture that this property holds for any specified sizes of the parts.

\begin{conj}
\label{general}
Let $G$ be a 2-connected graph of order $n$ and let $n = n_1+n_2+\cdots+n_k$ be any partition of
$n$. Then $V(G)$ can be partitioned into $k$ nearly connected parts $A_1, \ldots, A_k$ such
that $|A_i| = n_i$ for $1 \le i \le k$.
\end{conj}

If Conjecture \ref{general} is true, taking $n_i=r$ for all $i$, it implies that $G^r$ contains a 
$K_r$-factor for any 2-connected graph $G$ of order divisible by $r$. In Theorem~\ref{conn} we
prove a similar general statement for connected graphs. 

\begin{thm}
\label{conn}
Let $G$ be a connected graph of order $n$ and let $n = n_1+n_2+\cdots+n_k$ be a partition of $n$.
Then $V(G)$ can be partitioned into $k$ parts $A_1,A_2,\ldots,A_k$ such that $|A_i| = n_i$ and
there exists a subtree $T_i$ of $G$ of order at most $2n_i-1$ such that $A_i \subseteq V(T_i)$,
for $1 \le i \le k$.
\end{thm}

\noindent {\bf Proof:} It is sufficient to prove the result for trees. Let $T$ be a tree and fix
a vertex $u$ as the root. For any vertex $v \neq u$, let the subtree $T(v)$ of $T$ be the tree
obtained from $T$ by deleting the vertices in the component of $T-v$ that contains $u$, and let 
$T(u)$ be $T$ itself. Let $w$ be a vertex such that $|T(w)| \ge n_1$ and is farthest from 
$u$ in $T$. Then every component of $T(w)-w$ contains less than $n_1$ vertices. If $|T(w)| = n_1$,
let $A_1 = V(T(w))$. Since $T-A_1$ is a connected graph we can apply induction to find the
remaining parts. If $|T(w)| > n_1$, let $C_1,\ldots,C_k$ be a minimal set of components in 
$T(w)-w$ such that $|C_1|+\cdots+|C_k| \ge n_1$. Since $|C_i| < n_1$, we must have $|C_1|+\cdots
+|C_k| \le 2n_1-2$. Let $A_1$ be a set of $n_1$ vertices in $V(C_1) \cup \cdots \cup V(C_k)$ 
that are furthest from $w$ in $T(w)$. Then $A_1$ is a subset of vertices contained in the
subtree  $T_1$ of $T$ induced by $V(C_1) \cup \cdots \cup V(C_k) \cup \{w\}$ of order at most
$2n_1-1$, and $T-A_1$ is a connected subgraph. Again, the remaining parts can be found by
applying induction to $T-A_1$. \hfill $\Box$

It follows from Theorem~\ref{conn} that for any connected graph $G$ of order divisible by $r$,
$G^{2r-2}$ has a $K_r$-factor. The tree $T$ obtained from $K_{1,r+1}$ by subdividing each
edge $r-2$ times shows that in general $T^{2r-3}$ does not contain a $K_r$-factor. This tree
has $r^2$ vertices and $r+1$ leaf vertices at distance $2r-2$ from each other. Thus any 
partition of the vertex set into $r$ parts of size $r$ must have two leaf vertices at distance
$2r-2$ in the same part.

In Section~\ref{term} we introduce the terminology and notation needed for the proof of 
Theorem~\ref{main}. We assume all standard terms and define only those that are specific to this work.
Section~\ref{proof} gives the proof of Theorem~\ref{main}. We conclude in Section~\ref{conc} with 
some additional comments.

\section{Terminology}
\label{term}

In order to prove Theorem~\ref{main} by induction, we need to consider a larger class of graphs. 
We also need to consider partitions of a subset $A \subseteq V(G)$ of \emph{active} vertices
into nearly connected parts of size 4. The vertices not in $A$ are called \emph{dummy} vertices. 
These are needed in order to maintain the distances between the active vertices. 
Note that the fifth vertex in a subtree of order 5 that contains a nearly connected part
with four active vertices may or may not be a dummy vertex. Also, an active vertex may be
contained in many such subtrees, although it will be contained in exactly one nearly
connected part in the partition. 

We define certain subsets of rooted trees of small orders. Let $S_i$ for $0 \le i \le 3$ denote the set 
of all rooted trees with $i+1$ vertices. Let $S_i^+$ for $0 < i \le 3$ denote the set of rooted trees 
with $i+2$ vertices in which one of the vertices other than the root is marked as a dummy vertex. Let 
$S_2^-$ be the set containing only the path $P_3$ with the vertex of degree 2 as the root. Let $S_3^-$ 
be the set containing $P_4$ and $K_{1,3}$, both rooted at a vertex of degree more than 1. Let $S_5^-$ 
contain all rooted trees obtained from a tree in $S_3$ and a tree in $S_2$ by identifying their root nodes. 
If $G$ is a graph, $v$ a vertex in $V(G)$, and $S$ a set of rooted trees that are disjoint from $G$, let 
$G + S(v)$ denote the set of all graphs that can be obtained by identifying the root of some tree $T \in S$ 
with $v$. We say a graph in $G+S(v)$ is obtained by attaching the tree $T$ to $G$ at $v$. The vertices in 
$T$ other than the root are said to be attached to $v$. We define a partial order $\le $ on the sets of 
rooted trees, which is just the subset relation, except that we consider a tree in $S_i$ to also be in 
$S_i^+$.  Adding a dummy leaf vertex to a tree in $S_i$ gives a tree in $S_i^+$. 
Thus the relation $\le$ is defined by $S_1 \le S_1^+$, and $S_i^- \le S_i \le S_i^+$ for $i \in \{2,3\}$.

We consider non-trivial, connected multigraphs without a cutvertex called \emph{blocks}. The edges in 
the graph are assumed to be oriented arbitrarily and assigned labels according to the orientation. The 
label of an edge will be an element of the label set $\mathcal{L} = \{L_0, L_{00}, L_1, L_{10}, L_2, 
L_{20}, L_{21}, L_{30}, L_{31}, L_{32}\}$. We define an involution $f$ on the label set 
$\mathcal{L}$ such that if an edge $uv$ is labeled $L$ when oriented from $u$ to $v$, it is labeled $f(L)$ 
when oriented from $v$ to $u$. The involution $f$ is defined by $f(L) = L$ for all labels $L$ except 
$L_{31}$ and $L_{32}$, in which case, we define $f(L_{31}) = L_{32}$, and $f(L_{32}) = L_{31}$. 
An edge labeled $L_i$ or $L_{ij}$ is said to have weight $i$ for $0 \le i \le 3$. Each label $L$ is a set 
of ordered pairs $(P,Q)$ where $P,Q$ are sets of rooted trees and these are defined below. Note that 
the definitions are such that for any label $L$, $(P,Q) \in L$ if and only if $(Q,P) \in f(L)$.

\begin{enumerate}
\item
$L_0 = \{(S_0,S_0)\}$.
\item
$L_{00} = \{(S_0,S_0),(S_1,S_3^+),(S_1^+,S_3),(S_2,S_2^+),(S_2^+,S_2),(S_3,S_1^+),(S_3^+,S_1)\}$.
\item
$L_1 = \{(S_0,S_1), (S_1,S_0)\}$.
\item
$L_{10} = \{(S_0,S_1), (S_1,S_0), (S_2,S_3^+), (S_2^+,S_3), (S_3,S_2^+),(S_3^+,S_2) \}$.
\item
$L_2 = \{(S_0,S_2), (S_1,S_1), (S_2,S_0)\}$.
\item
$L_{20} = \{(S_0,S_2), (S_1,S_1), (S_2,S_0), (S_3,S_3^+), (S_3^+,S_3) \}$.
\item
$L_{21} = \{(S_0,S_2^-), (S_1,S_1^+), (S_1,S_5^-), (S_1^+,S_1), (S_5^-,S_1), (S_2^-,S_0), (S_3^-,S_3^-)\}$.
\item
$L_{30} = \{(S_0,S_3), (S_1,S_2), (S_2,S_1), (S_3,S_0)\}$.
\item
$L_{31} =\{(S_0,S_3^-), (S_1,S_2^-), (S_2,S_1^+), (S_2,S_5^-), (S_2^+,S_1), (S_3,S_0)\}$.
\item
$L_{32} =\{(S_0,S_3), (S_1,S_2^+), (S_1^+,S_2), (S_5^-,S_2), (S_2^-,S_1), (S_3^-,S_0)\}$.
\end{enumerate}

\begin{figure}
\includegraphics{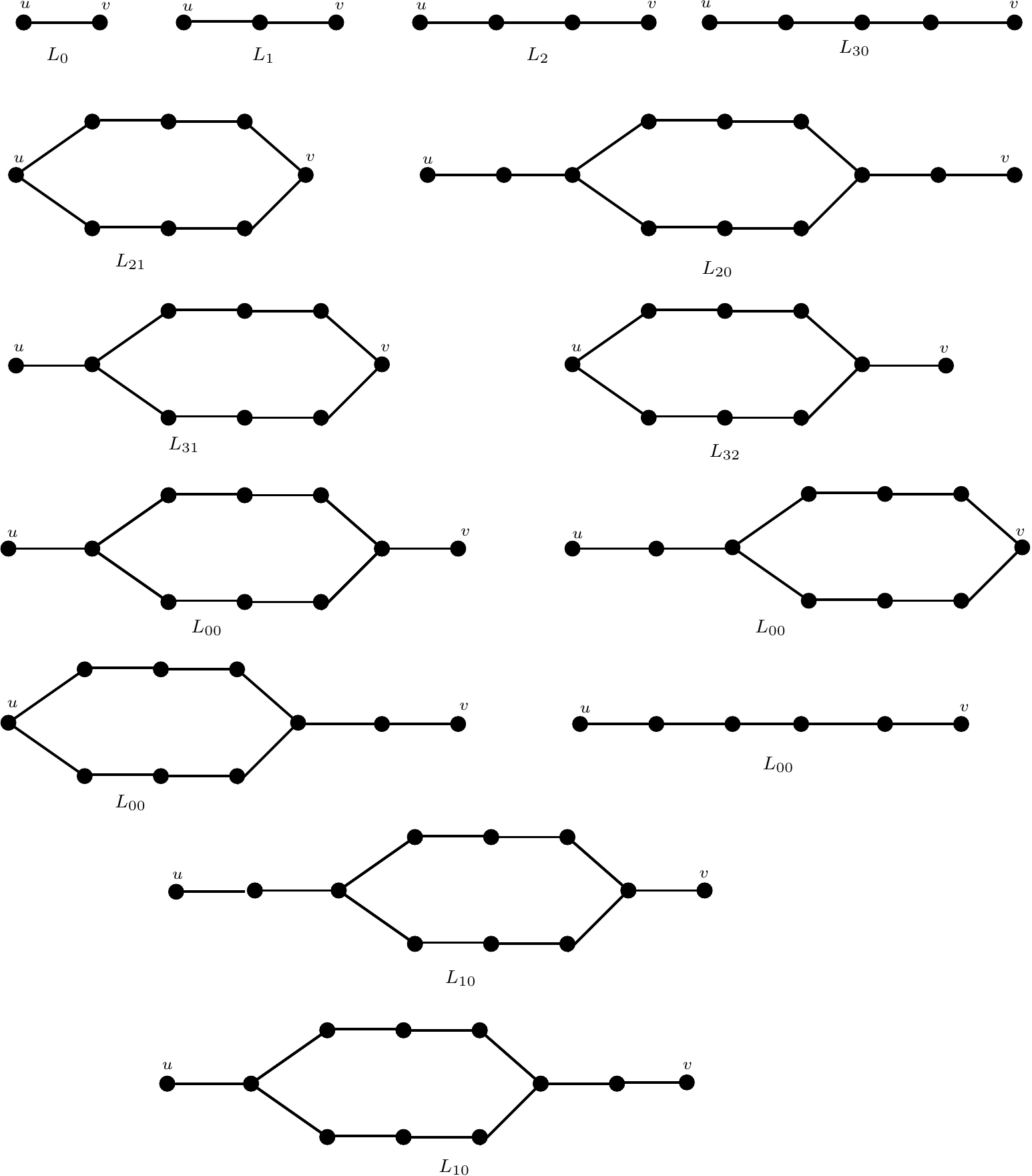}
\caption{Labels representing subgraphs}
\label{fig1}
\end{figure}

Informally, a labeled edge $e=uv$ represents a connected subgraph $H$ containing the vertices $u,v$
and the label of the edge indicates some properties of $H$. Suppose $e$ is oriented from
$u$ to $v$ and labeled $L$ with respect to this orientation. Then $(P,Q) \in L$ indicates that
there are disjoint subtrees $T_1 \in P$, $T_2 \in Q$ in $H$, rooted at $u$ and $v$ respectively,
such that the vertices in $(V(H) \setminus (V(T_1) \cup V(T_2))) \cup D$ can be partitioned into 
nearly connected parts of size 4, where $D$ is the set of dummy vertices in $T_1 \cup T_2$. 
Figure \ref{fig1} shows for each label an example of a subgraph that an edge with that label 
represents.

Given a graph $G$ with an orientation and assignment of labels to its edges, we define certain operations 
on the edges in $G$. Suppose $e=uv$ is an edge in $G$ oriented from $u$ to $v$ and having label $L$ with 
respect to this orientation. If $(P,Q) \in L$, any graph in $((G-uv)+P(u))+Q(v)$ is said to be obtained 
from $G$ by applying a $(P,Q)$ split to the edge $e$. Informally, this means that the edge $e$ is deleted 
from $G$ and some tree in $P$ ($Q$) is attached to $G$ at the vertex $u$ ($v$). For any edge with a label 
$L_{ij}$, the only operation allowed is applying a $(P,Q)$-split to the edge for some $(P,Q)  \in L_{ij}$. 
For an edge with label $L_i$, apart from applying a $(P,Q)$ split for some $(P,Q) \in L_i$, we also allow 
the operation of subdividing the edge exactly $i$ times. The new vertices added, if any, are considered to 
be active in the resulting graph. An edge $e$ with label $L(e)$ is said to admit a $(P,Q)$ split if there 
exist $P_1 \le P$ and $Q_1 \le Q$ such that $(P_1,Q_1) \in L(e)$.  We will also allow applying a $(P,Q)$ 
split to the edge $e$ if it admits it, although the actual split applied would be a $(P_1,Q_1)$ split.
Any graph obtained by a $(P_1,Q_1)$ split applied to $e$ can also be obtained by a $(P,Q)$ split, though 
the converse may not be true.

A split of $G$ is an assignment of an allowed operation to every edge in $G$. Given a split 
of $G$, we say a graph $H$ is obtained by splitting $G$ if it can be obtained by applying the specified 
operation to every edge in $G$. Note that different graphs can be obtained from the same split of $G$, 
depending on the trees attached to the vertices in $G$. In any such graph $H$, the vertices in $G$ are always 
considered to be active and the dummy vertices, if any, are in the trees attached to vertices in $G$. We say 
that $H$ is 4-partitionable if there exists a partition of the active vertices in $H$ into nearly connected 
parts of size 4. Suppose $v$ is an active vertex in such a graph $H$. Let $H'$ be obtained from $H$ by 
splitting the vertex $v$ into two vertices, an active vertex $v'$ and a dummy vertex $v''$, such that 
every vertex adjacent to $v$ in $H$ is adjacent to exactly one of $v'$ or $v''$ in $H'$. If $H'$ is 
4-partitionable, then so is $H$ and a 4-partition of $H$ is obtained from that of $H'$ by replacing the 
vertex $v'$ by $v$.    

A split of $G$ is said to be \emph{admissible} if any graph obtained by applying the specified
operation to every edge in $G$ is $4$-partitionable.

\section{Proof}
\label{proof}

We prove a more general statement that implies Theorem~\ref{main}. 

\begin{lem}
\label{split}
Let $G$ be a block of order $n$ whose edges are oriented and labeled. Let $w(G)$ be the sum of weights of
all edges in $G$ and suppose $w(G)+n$ is divisible by 4. Then $G$ has an admissible split.
\end{lem}

We first show that Lemma~\ref{split} implies Theorem~\ref{main}. Suppose $G$ is a 2-connected graph of 
order divisible by 4. Orient the edges in $G$ arbitrarily and assign the label $L_0$ to each edge. 
Lemma~\ref{split} implies that the resulting labeled graph has an admissible split. Since each edge is 
labeled $L_0$, the only operations allowed are applying an $(S_0,S_0)$ split, which is equivalent to 
deleting the edge, or subdividing the edge zero times, which is equivalent to retaining the edge. Thus 
there is a unique graph $H$ that can be obtained by applying a given split of $G$ and it is 
a spanning subgraph of $G$ with all vertices active. Since $H$ is $4$-partitionable by Lemma~\ref{split}, 
this implies $G$ itself is $4$-partitionable.

\noindent {\bf Proof of Lemma~\ref{split}:} The proof is by induction on the number of edges in $G$. 
Suppose $G$ has a single edge, which implies $G$ is $K_2$. Then the weight of the edge must be 2 and 
the edge must be labeled $L_2, L_{20}$ or $L_{21}$. If it is labeled $L_2$, then subdividing the 
edge twice gives $P_4$ with all vertices active and it is $4$-partitionable. If it is labeled $L_{20}$
or $L_{21}$, the edge admits an $(S_3,S_3^+)$ split. Applying this gives a graph with two components, 
one of which is a tree of order 4 and all vertices active, and the other is a tree of order 5 
with four active vertices. Any such graph is $4$-partitionable, hence the split is admissible.
We may now assume $G$ has at least two edges.

\noindent {\bf Case 1 Parallel Reduction.} 

Suppose $G$ contains a pair of multiple edges $e_1, e_2$ with endpoints $u,v$,  $w(e_1) = i$
and $w(e_2) = j$. Without loss of generality, we may assume $i \le j$ and both edges are oriented 
from $u$ to $v$. If $i = 0$, applying an $(S_0,S_0)$ split to $e_1$ gives the graph $G-e_1$ which 
satisfies the induction hypothesis, and has fewer edges. By induction, $G-e_1$ has an admissible 
split and hence so does $G$. 

We may assume $i > 0$. We construct a new labeled graph $G'$ by deleting the edges $e_1$ and 
$e_2$ and adding a new labeled edge $e$ with the same endpoints, oriented from $u$ to $v$. The 
label of $e$ is chosen such that $w(e) \equiv i+j \mbox{ mod } 4$. 
Then $G'$ is also a block with fewer edges that satisfies the induction hypothesis and has an 
admissible split. We show that an admissible split of $G$ can be obtained from an admissible split 
of $G'$. We apply the same operation in $G$ to all edges other than $e_1,e_2$ as in the admissible split 
of $G'$. The label of $e$ is chosen such that for any operation applied to $e$ in $G'$, we can choose 
appropriate allowed operations for $e_1$ and $e_2$, so that any graph obtained from $G$ by applying these
operations can also be obtained from $G'$ by applying the specified operation to $e$, possibly with some 
additional nearly connected parts. Since $G'$ has an admissible split, so does $G$.

\noindent {\bf Case 1.1} Both $e_1$ and $e_2$ are labeled $L_{30}$.

In this case let the label of $e$ be $L_{21}$. Suppose an $(S_0,S_2^-)$ split is applied to $e$. 
Apply an $(S_2,S_1)$ split to both $e_1$ and $e_2$. The vertices in the trees attached to $u$ form a 
nearly connected part and can be deleted to get a graph that can be obtained by an $(S_0,S_2^-)$
split applied to $e$. A symmetrical argument holds if an $(S_2^-,S_0)$ split is applied to $e$. Suppose 
an $(S_1,S_1^+)$ split is applied to $e$. Apply an $(S_1,S_2)$ split to $e_1$ and an $(S_0,S_3)$ split 
to $e_2$. Let $v'$ be a vertex adjacent to $v$ in the tree attached to $v$ that is obtained by splitting 
$e_1$. Deleting the vertices in the tree attached to $v$ obtained by splitting $e_2$, and considering 
$v'$ to be a dummy vertex, gives a graph that can be obtained by applying an $(S_1,S_1^+)$ split to $e$. 
Since this graph is 4-partitionable, adding the nearly connected part containing the deleted vertices and 
$v'$ gives a 4-partition of the graph obtained by splitting $G$. A symmetrical argument holds
if an $(S_1^+,S_1)$ split is applied to $e$. If an $(S_1,S_5^-)$ split is applied to $e$, apply
an $(S_1,S_2)$ split to $e_1$ and an $(S_0,S_3)$ split to $e_2$. A symmetrical argument holds
if an $(S_5^-,S_1)$ split is applied to $e$. If an $(S_3^-,S_3^-)$ split is applied to
$e$, apply an $(S_2,S_1)$ split to $e_1$ and an $(S_1,S_2)$ split to $e_2$ to get an
admissible split of $G$.

\noindent {\bf Case 1.2} Both $e_1$ and $e_2$ are labeled $L_1$.

In this case, let the label of $e$ be $L_2$. If in an admissible split of $G'$ an $(S_0,S_2)$ split 
is applied to $e$, apply an $(S_0,S_1)$ split to both $e_1$ and $e_2$.  This gives a graph that can
be obtained by an $(S_0,S_2)$ split applied to $e$ in $G'$. A symmetrical argument holds if an $(S_2,S_0)$ 
split is applied to $e$. If an $(S_1,S_1)$ split is applied to $e$, applying an $(S_1,S_0)$ split to 
$e_1$ and an $(S_0,S_1)$ split to $e_2$ gives an admissible split of $G$. The only remaining case 
is if the edge $e$ is subdivided twice in the split of $G'$. Then any graph $H'$ obtained by 
splitting $G'$ will contain a path $u,v_1,v_2,v$ where $v_1,v_2$ are the new vertices added. In 
this case, subdivide both edges $e_1$ and $e_2$ once in $G$. The resulting graph $H$ obtained 
by splitting $G$ will contain paths $u,v_1,v$ and $u,v_2,v$, instead of the path $u,v_1,v_2,v$. 
Any connected part of size 4 or 5 in $H'$ is also connected in $H$, since $v_1v_2$ is the 
only edge in $H'$ that is not in $H$, and any connected part of size 4 or 5 that contains $v_1$ and 
$v_2$ must contain either $u$ or $v$.  Since $H'$ is 4-partitionable, $H$ is also 4-partitionable and 
$G$ has an admissible split.

\noindent {\bf Case 1.3} All other cases.

We may assume, without loss of generality, that $e_1$ is not labeled $L_{30}$.
Let the label of $e$ be $L_{w(e)0}$. This would imply that the splits applied to $e$ are 
$(S_x,S_y)$ for some $x+y = w(e)$, or $(S_x,S_y^+)$ or $(S_x^+,S_y)$ for $x+y = w(e)+4$ for
$0 \le x,y \le 3$.

First suppose an $(S_x,S_y)$ split is applied to $e$ for some $x+y = i+j \le 3$, which implies $i=1$ and
$j \le 2$. Then either $x \le i \le j$ or $y \le j$. Without loss of generality we can assume $x \le j$, 
otherwise interchange $x$ and $y$ and the vertices $u$ and $v$. If $e_2$ admits an $(S_x,S_{j-x})$ split, 
apply this split to $e_2$. Applying an $(S_0,S_i)$ split to $e_1$ gives an admissible split of $G$. If this 
is not possible, then $e_2$ must be labeled $L_{21}$, $x=1$ and $y=2$. Then apply an $(S_1,S_0)$ split
to $e_1$ and an $(S_0,S_2^-)$ split to $e_2$ to get an admissible split of $G$.  

The same argument holds if an $(S_x,S_y^+)$ or $(S_x^+,S_y)$ split is applied to $e$, 
for some $x+y = i+j \ge 4$. Again, assuming $x \le j$, apply an $(S_x,S_{j-x}^+)$ split to $e_2$ and
an $(S_0,S_i)$ split to $e_1$ to get a graph that can be obtained by an $(S_x,S_y^+)$ split applied to 
$e$. A symmetrical argument can be used if an $(S_x^+,S_y)$ split is applied to $e$.

Suppose an $(S_x,S_y)$ split is applied to $e$ for $x+y = i+j-4 \le 2$. Suppose $x = 0$.
If $e_1$ admits an $(S_{i-y},S_y)$ split, apply it to $e_1$ and an $(S_j,S_0)$ split to $e_2$. The vertices 
attached to $u$ form a nearly connected part and can be deleted to get the required split. A symmetrical 
argument holds if $e_2$ admits an $(S_{j-y},S_y)$ split. If this is not possible, $e_1$ must be labeled 
$L_{21}$ or $L_{31}$ with $y = 1$ or $L_{32}$ with $y = 2$ and the same holds for $e_2$. 

If $y=1$ then $i+j = 5$ and $e_1$ must be labeled $L_{21}$ and $e_2$ must be labeled $L_{31}$.
Apply an $(S_0,S_2^-)$ split to $e_1$ and an $(S_0,S_3^-)$ split to $e_2$. One of the
vertices adjacent to $v$ obtained by splitting $e_1$ can be deleted along with the
three vertices attached to $v$ obtained by splitting $e_2$, to get a graph that can be
obtained by an $(S_0,S_1)$ split applied to $e$. If $y =2$, both $e_1$ and $e_2$ must
be labeled $L_{32}$. Apply an $(S_2^-,S_1)$ split to both $e_1$ and $e_2$. The vertices
attached to $u$ can then be deleted to get a graph that can be obtained by an $(S_0,S_2)$
split of $e$. A symmetrical argument holds if $y=0$.

Suppose $x=y=1$, which implies $i = j = 3$. If $e_1$ is labeled $L_{31}$, apply an $(S_1,S_2^-)$ 
split to $e_1$ and an $(S_0,S_3)$ split to $e_2$. Again, the vertices in the tree attached to $v$ 
obtained by splitting $e_2$, together with one of the vertices adjacent to $v$ obtained by splitting $e_1$, 
form a nearly connected part that can be deleted. This leaves a graph that can be obtained by an $(S_1,S_1)$ 
split applied to $e$ in $G'$. A symmetrical argument holds if $e_1$ is labeled $L_{32}$, in which case an 
$(S_2^-,S_1)$ split is applied to $e_1$ and an $(S_3,S_0)$ split to $e_2$.

The final case to consider is if $i=j=1$ and an $(S_3,S_3^+)$ or $(S_3^+,S_3)$ split
is applied to $e$. Since either $e_1$ or $e_2$, without loss of generality $e_1$, is not labeled
$L_1$, it admits an $(S_3,S_2^+)$ split and apply this split to $e_1$. Applying an $(S_0,S_1)$ split 
to $e_2$ gives graphs that can be obtained by an $(S_3,S_3^+)$ split applied to $e$. A symmetrical
argument holds if an $(S_3^+,S_3)$ split is applied to $e$. 

This completes all cases when $G$ has multiple edges. We may now assume $G$ is a simple graph.

\noindent {\bf Case 2 Series Reduction.}

Suppose $G$ has a vertex $v$ of degree 2. Let $e_1 = v_1v$ and $e_2 = vv_2$ be the two
edges incident with $v$ and let $w(e_1)=i$, $w(e_2)=j$.  Again assume $i \le j$, $e_1$ is oriented 
from $v_1$ to $v$ and $e_2$ from $v$ to $v_2$. Let $G'$ be the graph obtained from $G$ by deleting 
the vertex $v$ and adding an edge $e$ oriented from $v_1$ to $v_2$. We choose the label of $e$ so 
that $w(e)\equiv  i+j+1 \mbox{ mod } 4$. Then $G'$ is a block with fewer edges than $G$ and has an 
admissible split by induction. Again, the label of $e$ is chosen so that for any allowed operation 
applied to $e$ in $G'$, appropriate operations can be applied to $e_1$ and $e_2$ in $G$ to get an 
admissible split of $G$. If both the edges $e_1$ and $e_2$ are split, this will leave some trees
attached to the vertex $v$. In such cases, we ensure that the active vertices in these trees, along
with $v$, can be partitioned into nearly connected parts and deleting these parts gives
a graph that can be obtained from $G'$ by an admissible split. This gives an admissible  
split of $G$. In most cases, we only indicate the label of $e$ and the splits to be applied 
to $e_1,e_2$, for each possible split of $e$. Verifying that these are admissible splits of $G$ 
is a routine exercise.

\noindent {\bf Case 2.1} $e_1$ is labeled $L_0$ and $e_2$ is labeled $L_j$ for $j \in \{0,1\}$.

Let the label of $e$ be $L_{j+1}$. Suppose an $(S_x,S_y)$ split is applied to $e$ for some $0 \le 
x,y \le j+1$ and $x+y=j+1$. If $x=0$ then apply an $(S_0,S_0)$ split to $e_1$ and subdivide
the edge $e_2$ $j$ times. If $x > 0$, subdivide $e_1$ $0$ times and apply an $(S_{x-1},S_{j+1-x})$ 
split to $e_2$. If $e$ is subdivided $j+1$ times, then subdivide $e_1$ $0$ times and $e_2$ $j$ times.
It is easy to verify that in each case the operations applied to $e_1$ and $e_2$ give the same
graph as obtained by applying the operation to $e$.

\noindent {\bf Case 2.2} $e_1$ is labeled $L_0$ and $e_2$ is labeled $L_{21}$.

Let the label of $e$ be $L_{31}$. If an $(S_0,S_3^-)$ split is applied to $e$, apply an
$(S_0,S_0)$ split to $e_1$ and an $(S_3^-,S_3^-)$ split to $e_2$. For any other operation 
applied to $e$, the edge $e_1$ is subdivided $0$ times. The operation to be applied to $e_2$, 
depending on the operation applied to $e$, is specified in the table below. In each case, 
it is easy to verify that the resulting graph can be obtained by applying the specified operation 
to $e$.

$$
\begin{array}{|c|c|}
\hline
e & e_2 \\
\hline
(S_1,S_2^-) & (S_0,S_2^-) \\
\hline
(S_2,S_1^+) & (S_1,S_1^+) \\
\hline
(S_2,S_5^-) & (S_1,S_5^-) \\
\hline
(S_2^+,S_1) & (S_1^+,S_1) \\
\hline
(S_3,S_0) & (S_2^-,S_0)\\
\hline
\end{array}
$$

\noindent {\bf Case 2.3} $e_1$ is labeled $L_{00}$ and $e_2$ is labeled $L_{21}$.

Let the label of $e$ be $L_{31}$. In this case, both edges $e_1$ and $e_2$ will
be split and the operations to be applied are specified below. 

$$
\begin{array}{|c|c|c|}
\hline
e & e_1 & e_2 \\
\hline
(S_0,S_3^-) & (S_0,S_0) & (S_3^-,S_3^-) \\
\hline
(S_1,S_2^-) & (S_1,S_3^+) & (S_0,S_2^-) \\
\hline
(S_2,S_1^+) & (S_2,S_2^+) & (S_1,S_1^+) \\
\hline
(S_2,S_5^-) & (S_2,S_2^+) & (S_1,S_5^-) \\
\hline
(S_2^+,S_1) & (S_2^+,S_2) & (S_1^+,S_1) \\
\hline
(S_3,S_0) & (S_3,S_1^+) & (S_2^-,S_0) \\
\hline
\end{array}
$$

\noindent {\bf Case 2.4} $e_1$ is labeled $L_{21}$ and $e_2$ is labeled $L_{31}$.

In this case, label $e$ as $L_{21}$.  

$$
\begin{array}{|c|c|c|}
\hline
e & e_1 & e_2 \\
\hline
(S_0,S_2^-) & (S_0,S_2^-) & (S_1,S_2^-) \\
\hline
(S_1,S_1^+) & (S_1,S_1^+) & (S_2,S_1^+) \\
\hline
(S_1,S_5^-) & (S_1,S_1^+) & (S_2,S_5^-) \\
\hline
(S_1^+,S_1) & (S_1^+,S_1) & (S_2^+,S_1) \\
\hline
(S_2^-,S_0) & (S_2^-,S_0) & (S_3,S_0) \\
\hline
(S_3^-,S_3^-) & (S_3^-,S_3^-) & (S_0,S_3^-) \\
\hline
\end{array}
$$

\noindent {\bf Case 2.5} $e_1$ and $e_2$ are both labeled $L_{32}$.

In this case, label $e$ as $L_{32}$. 

$$
\begin{array}{|c|c|c|}
\hline
e & e_1 & e_2\\
\hline
(S_0,S_3) & (S_0,S_3) & (S_0,S_3) \\
\hline
(S_1,S_2^+) & (S_1,S_2^+) & (S_1,S_2^+) \\
\hline
(S_1^+,S_2) & (S_1^+,S_2) & (S_1^+,S_2) \\
\hline
(S_5^-,S_2) & (S_5^-,S_2) & (S_1^+,S_2) \\
\hline
(S_2^-,S_1) & (S_2^-,S_1) & (S_2^-,S_1) \\
\hline
(S_3^-,S_0) & (S_3^-,S_0) & (S_3^-,S_0) \\
\hline
\end{array}
$$

\noindent {\bf Case 2.6} All other cases.

In all other cases, let the label of $e$ be $L_{w(e)0}$.

\noindent {\bf Case 2.6.1} $i+j+1 \le 3$.

In this case, the possible splits that can be applied to $e$ are an $(S_x,S_y)$ split for
some $x+y=i+j+1 \le 3$, or an $(S_x,S_y^+)$ or $(S_x^+,S_y)$ split for $x+y = i+j+5 \le 6$.

Suppose an $(S_x,S_y)$ split is applied to $e$ for some $0 \le x,y \le 3$ and 
$x+y = i+j+1 \le 3$. Suppose $y \le j$, which implies $x \ge i+1$. If $e_2$ admits an $(S_{j-y},S_y)$ 
split, apply it to $e_2$. If $e_1$ is labeled $L_i$, subdivide $e_1$ $i$ times, otherwise apply an
$(S_x,S_{i+4-x}^+)$ split to $e_1$. Any edge of weight $i$ always admits such a split, unless it is 
labeled $L_i$. Again, if both edges are split, the active vertices in the trees attached to $v$, together 
with $v$, form a nearly connected part that can be deleted to get a graph that can be obtained by an 
$(S_x,S_y)$ split applied to $e$. 

The only cases where this is not possible is if $e_2$ is labeled $L_{21}$ or $L_{31}$ with $y=1$, 
or if $e_2$ is labeled $L_{32}$ with $y = 2$. Since $i+j+1 \le 3$ and  $j \ge 2$, we must have $i = 0$ 
and $j=2$. This implies $e_1$ must be labeled $L_0$ or $L_{00}$ and $e_2$ must be labeled $L_{21}$. These 
cases are considered in Cases 2.2 and 2.3, respectively. 

If $y > j$, then $x \le i$. Since $i+j+1 \le 3$ and $i \le j$, we have $i \le 1$. Then $e_1$
admits an $(S_x,S_{i-x})$ split and apply it to $e_1$. If $e_2$ is labeled $L_j$ subdivide it $j$ times
otherwise it admits an $(S_{j+4-y}^+,S_y)$ split, and in either case, this gives an admissible split of 
$G$. 

Suppose $i+j \le 1$ and an $(S_x,S_y^+)$ split is applied to $e$ for some $0 \le x,y \le 3$ and
$x+y = i+j+5$.  Then $i = 0$, $j \le 1$ and $x,y \ge 2$. If $e_1$ is labeled $L_i$ and $e_2$ is 
labeled $L_j$, then this case is considered in Case 2.1. Suppose $e_1$ is not labeled $L_0$. 
Then $e_1$ admits an $(S_x,S_{4-x}^+)$ split and apply this split to $e_1$. If $e_2$ is labeled
$L_j$ then subdivide $e_2$ $j$ times to get an admissible split of $G$. Otherwise $e_2$ admits
an $(S_{j+4-y},S_y^+)$ split, and applying  this gives an admissible split of $G$. If $e_1$ is
labeled $L_0$ and $e_2$ is not labeled $L_j$,  subdivide $e_1$ $0$ times and apply an $(S_{x-1},S_y^+)$ 
split to $e_2$ to get an admissible split of $G$. A symmetrical argument holds if an $(S_x^+,S_y)$ split 
is applied to $e$.

\noindent {\bf Case 2.6.2} $i+j+1 \ge 4$.

In this case $w(e) = i+j-3$ and the possible splits applied to $e$ are an $(S_x,S_y)$ split
for some $0 \le x,y \le 3$ and $x+y=i+j-3$ or an $(S_x,S_y^+)$ or $(S_x^+,S_y)$ split for
some $0 \le x,y \le 3$ and $x+y = i+j+1 \le 6$.

Suppose an $(S_x,S_y)$ split is applied to $e$ for some $x+y = i+j-3$. Then we have 
$x \le i$ and $y \le j$. If $e_1$ admits an $(S_x,S_{i-x})$ split then apply this split to $e_1$. 
The edge $e_2$ admits an $(S_{j-y}^+,S_y)$ split and applying it gives an admissible split
of $G$. A symmetrical argument holds if $e_2$ admits an $(S_{j-y},S_y)$ split. If neither of 
these is possible then $e_1$ and $e_2$ must be labeled one of $\{L_{21},L_{31},L_{32}\}$. If 
$e_1$ is labeled $L_{21}$ or $L_{32}$ then $x = 1$ and if it is labeled $L_{31}$ then $x = 2$. 
Similarly, if $e_2$ is labeled $L_{21}$ or $L_{31}$ (oriented from $v$ to $v_2$) then $y = 1$ 
and if it is labeled $L_{32}$ then $y = 2$. Since $x+y=i+j-3$, and $i \le j$, the possible 
combinations of $(L(e_1),L(e_2))$ are $(L_{21},L_{31})$, $(L_{32},L_{32})$ and $(L_{31},L_{31})$.
The third case is symmetric to the second, obtained by relabeling the edges $e_1$ and $e_2$.
The first case is considered in Case 2.4 and the second in Case 2.5.  

Suppose an $(S_x,S_y^+)$ or an $(S_x^+,S_y)$ split is applied to $e$, for some $x+y=i+j+1$ and 
$4 \le i+j+1 \le 6$. Then either $x \le i$ or $y \le j$. If $x \le i$, apply an $(S_x, S_{i-x}^+)$  
split to $e_1$.  If $e_2$ is labeled $L_j$ then subdivide it $j$ times, otherwise apply an 
$(S_{j+4-y},S_y^+)$ split to $e_2$. In either case, this gives an admissible split of $G$. If 
$y \le j$, apply an $(S_{j-y},S_y^+)$ split to $e_2$. If $e_1$ is labeled $L_i$, subdivide
it $i$ times, otherwise apply an $(S_x,S_{i+4-x}^+)$ split to $e_1$. In either case, this gives
an admissible split of $G$. A symmetrical argument holds if an $(S_x^+,S_y)$ split is applied to $e$.

This completes all cases when $G$ has a vertex of degree 2. We may now assume that $G$ is a 
simple 2-connected graph of minimum degree at least 3.

\noindent {\bf Case 3 Reducible edge.}

Let $e_1=vv_1$ be an edge in $G$, oriented from $v$ to $v_1$, with label $L_{32}$, such that
$G-e_1$ is 2-connected. We say $e_1$ is \emph{reducible} if there exists an edge $e_2=vv_2 \neq e_1$ 
in $G$, oriented from $v$ to $v_2$, with label $L_{32}$ or $L_{30}$. Suppose $G$ contains such an 
edge. Let $G'$ be the graph obtained from $G-e_1$ by relabeling the edge $e_2$ as $L_{20}$. Then 
$G'$ satisfies the induction hypothesis and has an admissible split. We show that an admissible 
split of $G$ can be obtained from that of $G'$. An $(S_3^-,S_0)$ split is applied to the edge 
$e_1$ in all cases. This results in a vertex adjacent to $v$ and a path of length 2 attached to 
$v$.

Suppose an $(S_0,S_2)$ split was applied to $e_2$ in $G'$. If $e_2$ is labeled $L_{30}$, apply
an $(S_1,S_2)$ split to $e_2$ in $G$. The vertices attached to $v$ resulting from the splits of 
$e_1$ and $e_2$ form a nearly connected part, and can be deleted. If $e_2$ is labeled $L_{32}$, 
apply an $(S_5^-,S_2)$ split to $e_2$ in $G$. The vertices attached to $v$ can now be partitioned 
into two nearly connected parts and can be deleted. In both cases, we get an admissible split of $G$. 
If an $(S_1,S_1)$ split is applied to $e_2$ in $G'$, apply an $(S_2,S_1)$ split to $e_2$ in $G$. 
Again, a nearly connected part containing four of the vertices attached to $v$ can be deleted,
leaving a single vertex adjacent to $v$. This gives a graph that can be obtained by
applying an $(S_1,S_1)$ split to $e_2$ in $G'$. If an $(S_2,S_0)$ split is applied to
$e_2$ in $G'$, apply an $(S_3,S_0)$ split to $e_2$ in $G$. Again deleting the vertices attached
to $v$ resulting from this split and the vertex adjacent to $v$ resulting from the split of
$e_1$, gives a graph that can be obtained by an $(S_2,S_0)$ split applied to $e_2$ in $G'$.
Finally, if an $(S_3^+,S_3)$ or $(S_3,S_3^+)$ split is applied to $e_2$ in $G'$, apply an
$(S_0,S_3)$ split to $e_2$ in $G$, to get a admissible split of $G$. This completes all cases, 
and we may assume $G$ has no reducible edge.

\noindent {\bf Case 4 Reducible vertex.}

Let $v$ be a vertex in $G$ of degree at least 3 such that $G-v$ is 2-connected. We say $v$
is \emph{reducible} if there exists at most one edge incident with $v$ that is labeled $L_{31}$ 
when oriented away from $v$.  Suppose there exists a reducible vertex $v$ in $G$.

If any edge $e$ incident with $v$ has weight $0$, apply an $(S_0,S_0)$ split to $e$ to get the graph 
$G-e$. Since $G-e$ satisfies the induction hypothesis, it has an admissible split, which gives an 
admissible split of $G$. We may assume every edge incident with $v$ has weight at least 1.

\noindent {\bf Case 4.1} Degree of $v$ is 3.

Suppose the degree of $v$ is exactly 3. Let $e_1=vv_1, e_2=vv_2, e_3=vv_3$ be the three edges incident
with $v$ and oriented away from $v$. Let the weights of $e_1,e_2,e_3$ be $i,j,k$ respectively, and assume 
that $i \ge j \ge k \ge 1$. If $i+j+k+1 = 4$, which implies $i=j=k=1$, apply an $(S_1,S_0)$ split to 
each of the three edges. The vertices attached to $v$, together with $v$, form a nearly connected part 
and can be deleted. The resulting graph is $G-v$, which satisfies the induction hypothesis and has an 
admissible split. This gives an admissible split of $G$. 

Suppose $i+j+k+1=8$. If $j = 2$ then we must have $i=3$ and $k=2$. Apply an $(S_2,S_0)$ split to $e_2$ 
and $e_3$. The vertices attached to $v$ obtained from these splits form a nearly connected part that 
can be deleted. Apply an $(S_3,S_0)$ split to $e_1$. The vertices attached to $v$ obtained by splitting 
$e_1$, along with $v$, form a nearly connected part that can be deleted. This gives the graph $G-v$ 
which has an admissible split by induction. The other possibility is that $i=j=3$ and $k=1$. In this case, 
apply an $(S_1,S_0)$ split to $e_3$, and an $(S_3,S_0)$ split to $e_1$ and $e_2$. The vertices attached to 
$v$, along with $v$, can be partitioned into two nearly connected parts, which can be deleted.  Again, the 
remaining graph is $G-v$ which has an admissible split by induction. This gives an admissible split of $G$.

\noindent {\bf Case 4.1.1} $5 \le i+j+k+1 \le 7$. 

In this case, we must have $i \ge 2$, $j \le 2$. 

\noindent {\bf Case 4.1.1.1} $e_1$ is labeled $L_{21}$ or $L_{32}$, $e_2$ is labeled $L_{21}$
and $k=1$.

Let $G'$ be the graph obtained from $G-v$ by adding an edge $e$ labeled $L_{i0}$, oriented from 
$v_1$ to $v_3$. Apply an $(S_2^-,S_0)$ split to $e_2$ in all cases. Suppose an $(S_x,S_y)$ split 
is applied to $e$ for some $x+y = i$. If $y \le 1$, apply an $(S_{1-y},S_y)$ split to $e_3$ and 
an $(S_{i-x}^+,S_x)$ split to $e_1$. The vertices attached to $v$, along with $v$, form a nearly 
connected part that can be deleted. The remaining graph can be obtained by an $(S_x,S_y)$ split 
applied to $e$. If $y = 2$, then $x = i-2$. Since $e_1$ is labeled either $L_{21}$ or $L_{32}$, 
in either case it admits an $(S_2^-,S_x)$  split. The vertices attached to $v$ obtained from the 
splits of $e_1$ and $e_2$ form a nearly connected part that can be deleted. If $e_3$ is labeled 
$L_1$, subdividing it once gives the required split, and if it is labeled $L_{10}$, apply an 
$(S_3^+,S_2)$ split to $e_3$. Again, the vertices attached to $v$ obtained from this split, 
along with $v$, form a nearly connected part that can be deleted.  If $y=3$ then $i = 3$,
$x=0$ and $e_1$ must be labeled $L_{32}$. Apply an $(S_3^-,S_0)$ split to $e_1$. The three vertices 
attached to $v$ resulting from this split, together with one of the vertices adjacent to $v$ obtained 
by splitting $e_2$, form a nearly connected part that can be deleted. Again, if $e_3$
is labeled $L_1$, subdivide it once, otherwise apply an $(S_2^+,S_3)$ split to
$e_3$ to get an admissible split in $G$.

\noindent {\bf Case 4.1.1.2} All other cases with $5 \le i+j+k+1 \le 7$.

Let $G'$ be the graph obtained from $G-v$ by adding an edge $e$ oriented from $v_1$ to $v_2$, with $w(e) =
i+j+k-3$ and label $L_{w(e)0}$. In all cases, an $(S_k,S_0)$ split is applied to $e_3$ in $G$. 

Suppose an $(S_x,S_y)$ split is applied to $e$ for some $x+y = i+j+k-3$. First consider the
case $x \le i$ and $y \le j$.   Suppose $e_1$ admits an $(S_{i-x},S_x)$ split. Apply this split to $e_1$ 
and an $(S_{j-y}^+,S_y)$ split to $e_2$. The vertices attached to $v$, along with $v$, form a nearly 
connected part that can be deleted. The remaining graph can be obtained by applying an $(S_x,S_y)$ 
split to $e$ in $G'$. The same argument holds if $e_2$ admits an $(S_{j-y},S_y)$ split. If this is not 
possible, $e_1$ must be labeled $L_{21}$ or $L_{31}$ with $x=1$ or $L_{32}$ with $x=2$. Similarly, 
$e_2$ must be labeled $L_{21}$ with $y=1$, since $j \le 2$.

Suppose $x=y=1$, which implies $i+j+k =5$, $k = 1$ and $i+j = 4$. Therefore both
$e_1$ and $e_2$ must be labeled $L_{21}$, $k=1$ and Case 4.1.1.1 can be applied. If $x=2$ and 
$y=1$, we have $i+j+k=6$, $i=3$, $j=2$, $k=1$, $e_1$ must be labeled $L_{32}$, $e_2$ is labeled $L_{21}$
and again Case 4.1.1.1 is applicable.

Suppose $x \le i$ but $y > j$. If $j = 2$, then $y=3$ which implies $x=0$ and
$i+k=4$. Apply an $(S_i,S_0)$ split to $e_1$. The vertices attached to $v$ obtained from this split
and those obtained by splitting $e_3$ form a nearly connected part that can be deleted. If $e_2$ is 
labeled $L_2$, subdivide $e_2$ two times, otherwise apply an $(S_3^+,S_3)$ split to $e_2$. 
This gives an admissible split of $G$. If $j=k=1$, then $x+y=i-1 \le 2$,
hence $x = 0$, $y= 2$ and $i=3$. Since $i+k = 4$, the same argument as when $j=2$ can be 
used. If $e_2$ is labeled $L_1$, subdivide it once, otherwise apply an $(S_3^+,S_2)$
split to get the required split of $G$. 

Finally if $x > i$, we must have $x = 3$, $i=j=k=2$ and $y=0$. Apply an $(S_2,S_0)$ split 
to both $e_2$ and $e_3$. The vertices attached to $v$ resulting from these splits form a 
nearly connected part that can be deleted. If $e_1$ is labeled $L_2$,
subdivide it twice to get the required split of $G$. Otherwise, apply an
$(S_3^+,S_3)$ split to $e_1$ and delete $v$ along with the vertices attached to it, to get
a graph that can be obtained from an $(S_3,S_0)$ split of $e$. 

The remaining possible splits of $e$ to be considered are an $(S_x,S_y^+)$ or $(S_x^+,S_y)$
split for some $x+y = i+j+k+1$ and $5 \le x+y \le 6$. In these cases, we must have $k = 1$
and $x,y \ge 2$ with at least one of them being 3. Apply an $(S_1,S_0)$ split to $e_3$.
Suppose $x=2$ and $y=3$, which implies $i=2$ and $j = 1$. Apply an $(S_0,S_2)$ split to $e_1$. If 
$e_2$ is labeled $L_1$, subdivide it once, otherwise apply an $(S_2^+,S_3)$ split to $e_2$ and 
delete the vertices attached to $v$ along with $v$. This gives a graph that can be obtained by an 
$(S_2,S_3)$ split of $e$, hence this argument can be used if either an $(S_2,S_3^+)$ or $(S_2^+,S_3)$ 
split is applied to $e$.

If $x=3$ and $y=2$, we again have $i=2$ and $j=1$. If $e_1$ admits an $(S_3^+,S_3)$ split, apply it 
to $e_1$. If $e_2$ is labeled $L_1$, subdivide it once, consider $v$ to be a dummy vertex, and 
delete the vertices attached to $v$ that are obtained by splitting $e_1$. The resulting graph can be 
obtained by an $(S_3,S_2^+)$ split applied to $e$ and is 4-partitionable. The deleted vertices along 
with $v$ form a nearly connected part, hence the given splits of $G$ are admissible.  If $e_2$ is 
not labeled $L_1$, it admits an $(S_3,S_2^+)$ split. Apply this split to $e_2$ and delete the three 
vertices attached to $v$ resulting from this split, along with the vertex adjacent to $v$ resulting 
from the split of $e_3$. The remaining vertices attached to $v$, resulting from the split of $e_1$, 
can be deleted along with $v$, to get a graph that can be obtained by an $(S_3,S_2^+)$ split of $e$. 
The other possibility is that $e_1$ is labeled $L_2$. In this case subdivide $e_1$ twice.  If $e_2$ 
admits an $(S_3,S_2^+)$ split, apply this split to $e_2$ and delete the vertices attached to $v$ 
resulting from this split, together with the vertex adjacent to $v$ obtained by splitting $e_3$. If 
$e_2$ is labeled $L_1$, subdivide it once and split the vertex $v$ into two, $v'$ and $v''$, where $v'$ 
is an active vertex adjacent to the vertex in the subdivision of $e_1$, and the vertex $v''$ is a 
dummy vertex adjacent to the vertex in the subdivision of $e_2$ and the vertex obtained by splitting 
$e_3$. This gives a graph $G'$ that can be obtained by an $(S_3,S_2^+)$ split of $e$, which is 
4-partitionable. In the 4-partition of $G'$, replace the vertex $v'$ by $v$ to get a 4-partition of $G$.

A similar argument holds if an $(S_3^+,S_2)$ split is applied to $e$. Apply an $(S_3,S_3^+)$ split to
$e_1$ if it is not labeled $L_2$, otherwise subdivide it twice. Apply an $(S_3^+,S_2)$ split to $e_2$ 
if it is not labeled $L_1$, otherwise subdivide it once. If both edges are split, the vertices
attached to $v$, along with $v$, can be partitioned into nearly connected parts and deleting these
gives a graph that can be obtained by an $(S_3^+,S_2)$ split applied to $e$. If only one of the edges $e_1,e_2$
is split, the vertices attached to $v$ obtained from this split can be deleted along with the vertex
obtained by splitting $e_3$. If $e_1$ is labeled $L_2$ and $e_2$ is labeled $L_1$, the vertex $v$ can be split 
into an active vertex $v'$ adjacent to the vertex in the subdivision of $e_2$, and a dummy vertex $v''$ 
adjacent to the other two vertices. This gives a graph that can be obtained by an $(S_3^+,S_2)$ split of $e$.

Finally, suppose $x=y=3$ and an $(S_3,S_3^+)$ split is applied to $e$. Again we have $k=1$ and apply
an $(S_1,S_0)$ split to $e_3$. We have either $i=3$ and $j=1$ or $i=j=2$. If $i = 3$, apply an
$(S_0,S_3)$ split to $e_1$. If $e_2$ is labeled $L_1$ subdivide it once, otherwise apply an 
$(S_2^+,S_3)$ split to $e_2$. This gives a graph obtained by an $(S_3,S_3)$ split applied to $e$,
and the same argument holds if an $(S_3^+,S_3)$ split is applied to $e$. If $i=j=2$, apply an $(S_3^+,S_3)$ 
split to $e_1$ if it is not labeled $L_2$, otherwise subdivide it twice. Apply an $(S_3,S_3^+)$ split
to $e_2$ if it is not labeled $L_2$, otherwise subdivide it twice, delete the vertices attached to $v$
obtained by splitting $e_1$, and mark $v$ as a dummy vertex. Again, unless both edges are labeled 
$L_2$, we can delete appropriate vertices attached to $v$ get a graph obtained by an $(S_3,S_3^+)$ split applied 
to $e$. If both $e_1$ and $e_2$ are labeled $L_2$, split the vertex $v$ into an active vertex $v'$ and a 
dummy vertex $v''$ such that $v'$ is adjacent to the vertex in the subdivision of $e_1$ while $v''$ is 
adjacent to the other two vertices. A symmetrical argument holds if an $(S_3^+,S_3)$ split is
applied to $e$.

This completes all cases when $5 \le i+j+k+1 \le 7$. 

\noindent {\bf Case 4.1.2} $i+j+k+1 = 9$. 

This implies $i=j=3$ and $k = 2$. Since there is at most one edge labeled $L_{31}$ incident with 
$v$ when oriented away from $v$, we may assume that if there is one, it is $e_1$, and thus $e_2$ 
is not labeled $L_{31}$.

\noindent {\bf Case 4.1.2.1} $e_1$ is labeled $L_{31}$, $e_2$ is labeled $L_{30}$ and $e_3$ is 
labeled $L_2$ or $L_{20}$. 

Let $G'$ be obtained from $G-v$ by adding an edge $e$ oriented from $v_2$ to $v_3$ labeled $L_{10}$. 
In this case an $(S_3,S_0)$ split is applied to $e_1$ in all cases. If an $(S_1,S_0)$ split is applied 
to $e$, apply a $(S_2,S_1)$ split to $e_2$, and an $(S_2,S_0)$ split to $e_3$. The vertices attached to 
$v$ obtained from the splits of $e_2$ and $e_3$ form a nearly connected part and can be deleted, and then 
the remaining vertices obtained from the split of $e_1$ can be deleted along with $v$. This gives the 
required split in $G$. If an $(S_0,S_1)$ split is applied to $e$, apply an $(S_3,S_0)$ split to $e_2$ 
and an $(S_1,S_1)$ split to $e_3$. The argument is now the same as in the previous case. If an 
$(S_2,S_3^+)$ split is applied to $e$, apply an $(S_1,S_2)$ split to $e_2$. The vertices attached to $v$
obtained from the split of $e_1$ can be deleted along with the vertex attached to $v$ obtained from the 
split of $e_2$. If $e_3$ is labeled $L_{20}$, apply an $(S_3^+,S_3)$ split to $e_3$ and delete the 
vertices attached to $v$, along with $v$. If $e_3$ is labeled $L_2$, subdivide the edge $e_3$ twice. 
In either case, the resulting graph can be obtained from an $(S_2,S_3)$ split applied to $e$, and hence 
the same argument holds if an $(S_2^+,S_3)$ split is applied to $e$. If an $(S_3,S_2^+)$ split is applied 
to $e$, apply an $(S_0,S_3)$ split to $e_2$, an $(S_0,S_2)$ split to $e_3$ and delete the vertices 
attached to $v$, along with $v$. The resulting graph can be obtained by an $(S_3,S_2)$ split applied to 
$e$, hence the same argument holds if an $(S_3^+,S_2)$ split is applied to $e$.

\noindent {\bf Case 4.1.2.2}  $e_1$ is labeled $L_{31}$, $e_2$ is labeled $L_{32}$ and $e_3$ is labeled 
$L_2$ or $L_{20}$. 

Let $G'$ be obtained from $G-v$ by adding an edge $e$ oriented from $v_1$ to $v_3$, with label 
$L_{10}$. An $(S_3^-,S_0)$ split is applied to $e_2$ which results in a vertex adjacent to $v$ and 
another path of length two attached to $v$. If an $(S_1,S_0)$ split is applied to $e$, apply an 
$(S_2^+,S_1)$ split to $e_1$ and an $(S_2,S_0)$ split to $e_3$. The two vertices in the path of length 
two attached to $v$ resulting from the split of $e_2$ form a nearly connected part with the two vertices
obtained from the split of $e_3$. These can be deleted and the remaining vertices attached to $v$ also form 
a nearly connected part, along with $v$, and can all be deleted. This gives a graph that can be obtained by 
an $(S_1,S_0)$ split applied to $e$. If an $(S_0,S_1)$ split is applied to $e$, apply an $(S_3,S_0)$ split 
to $e_1$ and an $(S_1,S_1)$ split to $e_3$. It can be checked that the vertices attached to $v$ can be 
partitioned into two nearly connected parts which can be deleted to get an admissible split of $G$. If an 
$(S_2,S_3^+)$ split is applied to $e$, apply an $(S_1,S_2^-)$ split to $e_1$. The vertices attached to $v$ 
obtained from the splits of $e_1$ and $e_2$ form a nearly connected part that can be deleted. If $e_3$ is 
labeled $L_{20}$ apply an $(S_3^+,S_3)$ split to $e_3$ and delete the vertices attached to $v$ along with 
$v$. This gives a graph that can be obtained by an $(S_2,S_3)$ split applied to $e$. If $e_3$ is labeled
$L_2$, subdivide it twice to get the same graphs. The same argument holds for an $(S_2^+,S_3)$ 
split of $e$. If an $(S_3,S_2^+)$ split is applied to $e$, apply an $(S_0,S_3^-)$ split to $e_1$, and
an $(S_0,S_2)$ split to $e_3$. Again, the vertices attached to $v$, along with $v$, can be deleted to
get a graph obtained by an $(S_3,S_2)$ split applied to $e$ and the same argument can be used 
if an $(S_3^+,S_2)$ split is applied to $e$. 

\noindent {\bf Case 4.1.2.3} All other cases with $i+j+k+1 = 9$.

Let $G'$ be the graph obtained from $G-v$ by adding an edge $e$ oriented from $v_1$ to $v_2$ with 
label $L_{10}$.  Suppose $e_1$ is also not labeled $L_{31}$. In all cases, apply an $(S_2,S_0)$ split 
to $e_3$. Suppose an $(S_1,S_0)$ split is applied to $e$. Since $e_1$ is not labeled $L_{31}$ it admits
an $(S_2,S_1)$ split. The vertices attached to $v$ obtained from the splits of $e_1$ and $e_3$ 
form a nearly connected part that can be deleted. Applying an $(S_3,S_0)$ split to $e_2$, and 
deleting the vertices attached to $v$, along with $v$, gives the required split of $G$. A 
symmetrical argument holds if an $(S_0,S_1)$ split is applied to $e$. If an $(S_2,S_3^+)$ split 
is applied to $e$, apply an $(S_1^+,S_2)$ split to $e_1$ and an $(S_0,S_3)$ split 
to $e_2$. The vertices attached to $v$ can now be deleted, along with $v$, to get a graph that can 
be obtained by an $(S_2,S_3)$ split of $e$. Thus the same argument holds for an $(S_2^+,S_3)$ split 
applied to $e$ and a symmetrical argument holds if an $(S_3,S_2^+)$ or $(S_3^+,S_2)$ split is 
applied to $e$.  

Suppose $e_1$ is labeled $L_{31}$.  Since $e_2$ is labeled either $L_{30}$ or $L_{32}$, we may
assume $e_3$ is labeled $L_{21}$, otherwise either Case 4.1.2.1 or 4.1.2.2 applies. Apply an 
$(S_2^-,S_0)$ split to $e_3$. In this case, the same argument as in the case when $e_1$ is
not labeled $L_{31}$ holds, except when an $(S_1,S_0)$ split is applied to $e$. In this case, an 
$(S_2^+,S_1)$ split can be applied to $e_1$ and an $(S_3,S_0)$ split to $e_2$. The vertices 
attached to $v$ resulting from the split of $e_2$ can be deleted along with one of the two 
vertices adjacent to $v$ that result from splitting $e_3$. The remaining vertices attached to $v$, 
along with $v$, can now be deleted to get the required split graph. For all other splits 
applied to $e$, the argument is the same as in the case when $e_1$ is not labeled $L_{31}$.

\noindent {\bf Case 4.1.3} $i+j+k+1=10$, which implies $i=j=k=3$.

We may assume that if there is any edge labeled $L_{31}$, it is $e_1$. Suppose at least one of 
$e_1,e_2,e_3$ is labeled $L_{32}$, and without loss of generality, assume it is $e_3$. Then since 
$e_2$ is labeled $L_{30}$ or $L_{32}$, $e_3$ is a reducible edge and we can apply Case 3. The 
only other possibility is that both $e_2,e_3$ are labeled $L_{30}$ and $e_1$ is labeled either 
$L_{30}$ or $L_{31}$. In this case let $G'$ be obtained from $G-v$ by adding an edge $e$ with 
label $L_{20}$ oriented from $v_2$ to $v_3$. The edge $e_1$ is split by applying an $(S_3,S_0)$ 
split.

The splits to be applied to $e_2$ and $e_3$ are shown below for each possible split of
$e$. In all cases, it can be verified that the vertices attached to $v$, along with $v$,
can be partitioned into nearly connected parts and deleted to a get graphs that can be 
obtained from the specified split of $e$.

$$
\begin{array}{|c|c|c|}
\hline
e & e_2 & e_3 \\
\hline
(S_0,S_2) & (S_3,S_0) & (S_1,S_2) \\
\hline
(S_1,S_1) & (S_2,S_1) & (S_2,S_1) \\
\hline
(S_2,S_0) & (S_1,S_2) & (S_3,S_0) \\
\hline
(S_3,S_3^+) & (S_0,S_3) & (S_0,S_3) \\
\hline
(S_3^+,S_3) & (S_0,S_3) & (S_0,S_3)\\
\hline
\end{array}
$$

This completes all cases when there exists a reducible vertex of degree 3.

\noindent {\bf Case 4.2} Degree of $v$ is at least 4.

Let the degree of $v$ be $d \ge 4$ and let $e_i = vv_i$ for $1 \le i \le d$, be the edges
incident with $v$, oriented from $v$ to $v_i$. As argued previously, we can assume $w(e_i) > 0$
for $1 \le i \le d$. Suppose there exist two edges $e_i,e_j$ such that $w(e_i)+w(e_j) = 4$.
Apply an $(S_{w(e_i)},S_0)$ split to $e_i$ and an $(S_{w(e_j)},S_0)$ split to $e_j$. The
vertices attached to $v$ obtained from these splits form a nearly connected part and can
be deleted. The remaining graph is $G-\{e_i,e_j\}$ which satisfies the induction hypothesis
and has an admissible split. This gives an admissible split of $G$.

Therefore there can be at most one edge of weight 2, and either there is
no edge of weight 1 or no edge of weight 3 incident with $v$. 

\noindent {\bf Case 4.2.1} There is no edge of weight 3 incident with $v$. 

Assume that if there is an edge of weight 2, it is $e_1$. If $d \ge 6$, apply an 
$(S_1,S_0)$ split to the edges $e_{d-3}, e_{d-2},e_{d-1}$ and $e_d$. The four vertices attached to 
$v$ resulting from these splits form a nearly connected part and can be deleted. The remaining
graph $G-\{e_{d-3}, e_{d-2},e_{d-1},e_d\}$ satisfies the induction hypothesis and has an
admissible split. This gives an admissible split of $G$. If $d=5$ and $e_1$ has weight 2, we
can use the same argument by applying an $(S_2,S_0)$ split to $e_1$ and an $(S_1,S_0)$ 
split to $e_4$ and $e_5$. The only remaining cases are if $d=5$ and all edges have weight 1
and if $d=4$. 

If $d=5$ or $d=4$ and $w(e_1) = 2$, let $G'$ be the graph obtained from $G-v$
by adding an edge with label $L_{20}$, oriented from $v_1$ to $v_2$. Apply an $(S_1,S_0)$
split to $e_3,e_4$ and if $d=5$, also to $e_5$. Suppose an $(S_2,S_0)$ split is applied to $e$. 
If $d=4$, apply an $(S_0,S_2)$ split to $e_1$ and an $(S_1,S_0)$ split to $e_2$. The vertices attached to
$v$, along with $v$, form a nearly connected part that can be deleted to get the
required split. If $d=5$ and $e_1$ is not labeled $L_1$, apply an $(S_3^+,S_2)$ split
to $e_1$, and an $(S_1,S_0)$ split to $e_2$. The four vertices adjacent to $v$,
resulting from the splits of $e_2,e_3,e_4,e_5$ form a nearly connected part and
can be deleted. The vertices attached to $v$ due to the split of $e_1$ form a
nearly connected part along with $v$ and can be deleted to get the required split.
If $e_1$ is labeled $L_1$, subdivide it once and use the same argument to get the required split. 
If an $(S_0,S_2)$ split is applied to $e$, apply an $(S_{w(e_1)},S_0)$ split to $e_1$, and delete 
the vertices attached to $v$ obtained by this, along with those obtained by splitting $e_3,e_4$ and
$e_5$ if $d=5$.  If $e_2$ is labeled $L_{10}$, apply an $(S_3^+,S_2)$ split to $e_2$, otherwise 
subdivide it once to get an admissible split of $G$.
Suppose an $(S_1,S_1)$ split is applied to $e$. If $d=5$ apply $(S_0,S_1)$ splits to 
both $e_1$ and $e_2$, and if $d=4$, apply an $(S_1^+,S_1)$ split to $e_1$ and an $(S_0,S_1)$ 
split to $e_2$. In both cases, the vertex $v$ along with the vertices attached to it forms a 
nearly connected part and can be deleted to get the required split of $G$. 
Suppose an $(S_3,S_3^+)$ split is applied to $e$. If $e_2$ is not labeled $L_1$, apply 
an $(S_2,S_3^+)$ split to $e_2$, otherwise subdivide it once. If $d=5$, apply an 
$(S_2^+,S_3)$ split to $e_1$ if it is not labeled $L_1$ otherwise subdivide it once. 
If $d=4$, apply an $(S_3^+,S_3)$ split to $e_1$ if it is not labeled $L_2$, otherwise 
subdivide it twice. It can be checked that unless both $e_1$ and $e_2$ are subdivided, 
appropriate nearly connected parts can be deleted to get the required split of $G$. 
In case both edges are subdivided, split the vertex $v$ into an active vertex $v'$ 
adjacent to the vertex in the subdivision of $e_1$ and a dummy vertex $v''$ adjacent to the 
vertex in the subdivision of $e_2$. Two of the remaining neighbors of $v$ are connected to $v''$
and if $d=5$, the remaining one is connected to $v'$. This gives the required split.  A symmetrical 
argument holds if an $(S_3^+,S_3)$ split is applied to $e$. In this case an $(S_2^+,S_3)$ split is applied 
to $e_2$ if it is not labeled $L_1$ else it is subdivided once. If $w(e_1) = 1$, either an $(S_2,S_3^+)$ 
split is applied to $e_1$ or it is subdivided once, and if $w(e_1) =2$, either an $(S_3,S_3^+)$ split 
is applied to $e_1$ or it is subdivided twice. In case both edges are subdivided, split the vertex $v$ into 
a dummy vertex $v'$ and an active vertex $v''$ as in the previous case. This gives a graph that can be 
obtained by an $(S_3^+,S_3)$ split applied to $e$. 

The remaining case is if $d=4$ and all edges have weight 1. Let $G'$ be obtained from $G-v$ by adding
an edge $e$ oriented from $v_1$ to $v_2$ with label $L_{10}$. If an $(S_0,S_1)$ split is applied
to $e$, apply an $(S_1,S_0)$ split to $e_1$ and an $(S_0,S_1)$ split to $e_2$. This gives an
admissible split of $G$. A symmetrical argument holds if an $(S_1,S_0)$ split is applied to $e$.
If an $(S_2,S_3^+)$ split is applied to $e$, apply an $(S_3^+,S_2)$ split to $e_1$ if it is not
labeled $L_1$ otherwise subdivide it once. Apply an $(S_2,S_3^+)$ split to $e_2$ if it is not
labeled $L_1$ otherwise subdivide it once. If both edges are subdivided, split the vertex $v$
into an active vertex $v'$ adjacent to the vertex in the subdivision of $e_1$ and a dummy
vertex $v''$ adjacent to all other neighbours of $v$. This gives a graph that can be obtained
by an $(S_2,S_3^+)$ split applied to $e$. Exactly symmetrical arguments hold if any of
$(S_2^+,S_3)$, $(S_3,S_2^+)$, or $(S_3^+,S_2)$ splits are applied to $e$.

\noindent {\bf Case 4.2.2} There are no edges of weight 1 incident with $v$. 

Again, if there is any edge with label $L_{31}$, let it be $e_1$ and if there is an edge of weight 2, 
let it be $e_2$. Then $e_3$ and $e_4$ have weight 3 and neither of them is labeled $L_{31}$. If either of
them is labeled $L_{32}$ it is a reducible edge in $G$, a contradiction. Therefore both $e_3$ and $e_4$ 
must be labeled $L_{30}$. Let $G'$ be the graph obtained from $G-\{e_3,e_4\}$ by adding an edge $e$ with 
label $L_{20}$ oriented from $v_3$ to $v_4$. Then $G'$ satisfies the induction hypothesis and has an 
admissible split. We show that an admissible split of $G$ can be obtained from that of $G'$. For each 
possible split of $e$, the splits to be applied to $e_3$ and $e_4$ are shown below. In all cases,
the vertices attached to $v$ resulting from these splits, if any, form a nearly connected part,
and can be deleted to get a graph that can be obtained from the corresponding split of
$e$. 

$$
\begin{array}{|c|c|c|}
\hline
e & e_3 & e_4 \\
\hline
(S_0,S_2) & (S_3,S_0) & (S_1,S_2) \\
\hline
(S_1,S_1) & (S_2,S_1) & (S_2,S_1) \\
\hline
(S_2,S_0) & (S_1,S_2) & (S_3,S_0) \\
\hline
(S_3,S_3^+) & (S_0,S_3) & (S_0,S_3) \\
\hline
(S_3^+,S_3) & (S_0,S_3) & (S_0,S_3) \\
\hline
\end{array}
$$

This completes the list of reducible configurations. We show that $G$ must contain one of these.  
We may assume $G$ is a simple 2-connected graph of minimum degree at least 3, otherwise either parallel
or series reduction can be applied to $G$.  If $G$ is 3-connected 
then $G-v$ is 2-connected for every vertex $v$. If there is no reducible vertex then for every vertex 
$v$, there must be at least two edges incident with $v$ that are labeled $L_{31}$ when oriented away 
from $v$. This implies there exists a vertex $v$ with two edges incident with $v$ that are labeled 
$L_{32}$ when oriented away from $v$. Then both these edges are reducible. 

If $G$ is not 3-connected, let $\{u,v\}$ be a 2-cut such that the order of the smallest
component $C$ of $G-\{u,v\}$ is as small as possible. Since the minimum degree of $G$ is
at least 3, $C$ must contain at least two vertices and for any vertex $x \in V(C)$, $G-x$
is 2-connected. Again, for every vertex $x \in V(C)$, there must be at least two edges incident
with $x$ that are labeled $L_{31}$ when oriented away from $x$. This implies that either
there is a vertex $y \in V(C)$ such that there are two edges in $C$ incident with $y$ that
are labeled $L_{32}$ when oriented away from $y$, or for some  $y \in \{u,v\}$,
there are two edges labeled $L_{32}$ when oriented from  $y$ to vertices in $C$.
In either case, we get a reducible edge in $G$.

This completes the proof of Lemma~\ref{split} and hence that of Theorem~\ref{main}.
\hfill $\Box$
\section{Remarks}
\label{conc}

We show that if true, Conjecture \ref{general} would be best possible as there exist
2-connected graphs $G$ of order divisible by $r$ such that $G^{r-1}$ does not contain
a $K_r$-factor.

Let $G$ be the graph obtained from $K_4$ by subdividing five of the edges $r-1$ times 
and one edge $r+1$ times. Then $G$ has $5(r-1)+r+1+4 = 6r$ vertices, and we claim that
$G^{r-1}$ does not contain a $K_r$-factor. Suppose $G^{r-1}$ can be partitioned into six
cliques of order $r$. Since the vertices of degree 3 are at distance $\ge r$ from each other in $G$,
they must be in different cliques, and there are two cliques that do not contain any of these
vertices. At most one of these can be contained in the $r+1$ vertices that subdivide one
edge of $K_4$.  Hence there exists such a clique that contains vertices from the
subdivision of more that one edge. The edges of $K_4$ whose subdivision vertices are contained
in this clique must have a common endpoint $v$ in $K_4$, otherwise the distance between them
is greater than $r$. The clique that contains the vertex $v$ must include vertices from
the subdivisions of the three edges in $K_4$ incident with $v$. All these vertices must be at
distance at most $r-1$ from $v$. However, if $2r-1$ vertices at distance at most $r-1$ from $v$
are partitioned into two parts, some part must contain two vertices at distance at least $r$
from each other.  

If $r$ is even, a simpler example is the graph consisting of $r+2$ internally disjoint 
paths of length $r$ between two vertices $u,v$. Then the centers of the 
$r+2$ paths are at distance $r$ from each other in $G$, and any partition of the $r^2+r$ vertices 
in $G$ into $(r+1)$ parts of size $r$ must contain two of them in the same part.

It would be interesting to consider what happens for graphs with higher connectivity. In
particular, let $f(k,r)$ be the smallest integer such that for any $k$-connected graph $G$
of order divisible by $r$, $G^{f(k,r)}$ contains a $K_r$-factor. While we do not know the
exact behavior of even $f(3,r)$, it was conjectured in~\cite{AD} that $f(3,3) = 2$. As far
as we know, even this is not resolved. It may be noted though that assuming a higher
edge-connectivity does not help, and for all $k \ge 1$ and sufficiently large $r$, there exist
$k$-edge connected graphs $G$ of order divisible by $r$ such that $G^{2r-3}$ does not contain
a $K_r$-factor.

\end{document}